\newcommand{\arx}[1]{\texttt{http://arxiv.org/abs/#1}}

\documentclass{elsart}
\usepackage{amssymb, latexsym}

\renewcommand{\(}{\left (}
\renewcommand{\)}{\right )}
\newcommand{\p}{\mathfrak p}
\newcommand{\R}{{\mathbb R}}
\newcommand{\dx}{{dx}}
\newcommand{\be}{\begin{enumerate}}
\newcommand{\ee}{\end{enumerate}}
\newcommand{\bi}{\begin{itemize}}
\newcommand{\ei}{\end{itemize}}
\long\def\forget#1\forgotten{}

\begin{document}
\begin{frontmatter}

\title{Bernoulli numbers and the probability of a birthday surprise\thanksref{lea}}
\thanks[lea]{Dedicated to my wife Lea on her birthday}

\author{Boaz Tsaban}
\address{Department of Mathematics and Computer Science, Bar-Ilan University,
Ramat-Gan 52900, Israel} \ead{tsaban@macs.biu.ac.il}
\ead[url]{http://www.cs.biu.ac.il/\~{}tsaban}

\begin{abstract}
A birthday surprise is the event that, given $k$ uniformly random
samples from a sample space of size $n$, at least two of them are identical.
We show that Bernoulli numbers can be used to derive arbitrarily exact
bounds on the
probability of a birthday surprise.
This result can be used in arbitrary precision calculators, and it can
be applied to better understand some questions in
communication security and pseudorandom number
generation.
\end{abstract}

\end{frontmatter}

\section{Introduction}

In this note we address the probability $\beta^k_n$ that in a sample of
$k$ uniformly random
elements out of a space of size $n$ there exist at least two identical
elements. This problem has
a long history and a wide range of applications.
The term \emph{birthday surprise} for a collision
of (at least) two elements in the sample comes from the case $n=365$,
where the problem can be
stated as follows: Assuming that the birthday of people distributes
uniformly over the year,
what is the probability that in a class of $k$ students, at least
two have the same birthday?

It is clear (and well known) that the expected number of collisions
(or birthdays) in a sample
of $k$ out of $n$ is:
$${k \choose 2}\frac{1}{n} = \frac{k(k-1)}{2n}.$$
(Indeed, for each distinct $i$ and $j$ in the range $\{1,\dots,k\}$, let
$X_{ij}$ be the random variable taking the value $1$ if samples
$i$ and $j$ obtained the same value and $0$ otherwise.
Then the expected number of collisions is
$E(\sum_{i\neq j} X_{ij})=\sum_{i\neq j}E(X_{ij})=\sum_{i\neq j}\frac{1}{n}=
{k \choose 2}\frac{1}{n}$.)

Thus,
$28$ students are enough to make the expected number of common birthdays
greater than $1$.
This seemingly surprising phenomenon has
got the name \emph{birthday surprise}, or \emph{birthday paradox}.

In several applications, it is desirable to have exact bounds on the
probability of a collision.
For example, if some electronic application chooses pseudorandom
numbers as passwords for its
users, it may be a \emph{bad} surprise if two users get the same
password by coincidence.
It is this term ``by coincidence'' that we wish to make precise.

\section{Bounding the probability of a birthday surprise}

When $k$ and $n$ are relatively small, it is a manner of simple
calculation to determine
$\beta^k_n$. The probability that all samples are distinct is:
\begin{equation}
\pi^k_n = \(1-\frac{1}{n}\)\cdot \(1-\frac{2}{n}\)\cdot \ldots \cdot
\(1-\frac{k-1}{n}\),
\nonumber
\end{equation}
and $\beta^k_n = 1-\pi^k_n$.
For example, one can check directly that $\beta^{23}_{365}>1/2$,
that is, in a class of $23$ students the probability that
two share the same birthday is greater than $1/2$. This is another
variant of the \emph{Birthday surprise}.%
\footnote{To experience this phenomenon experimentally, the reader is
referred to \cite{SHweb}.}

The calculation becomes problematic when $k$
and $n$ are large,
both due to precision problems and computational complexity
(in cryptographic applications $k$ may be of the order of trillions, i.e.,
thousands of billions).
This problem can be overcome by considering the logarithm of the product:
$$\ln(\pi^k_n) = \sum_{i=1}^{k-1} \ln\(1-\frac{i}{n}\).$$
Since each $i$ is smaller than $n$, we can use the Taylor expansion
$\ln(1-x)=-\sum_{m=1}^\infty x^m/m$ ($|x|<1$)
to get that
\begin{equation}\label{sum}
-\ln(\pi^k_n) = \sum_{i=1}^{k-1} \sum_{m=1}^\infty\frac{(i/n)^m}{m} =
\sum_{m=1}^\infty \frac{1}{mn^m} \sum_{i=1}^{k-1} i^m.
\end{equation}
(Changing the order of summation is possible because the sums involve
positive coefficients.)

The coefficients $\p(k-1,m):=\sum_{i=1}^{k-1} i^m$ (which are often
called \emph{sums of powers}, or
simply \emph{power sums})
play a key role in our estimation of the birthday probability.
Efficient calculations of the first few power sums go back to
ancient mathematics.%
\footnote{
Archimedes (ca.\ 287-212 BCE) provided a geometrical derivation of a
``formula'' for the
sum of squares \cite{VKATZ}.
}
In particular, we have:
$\p(k,1)= k(k+1)/2$, and $\p(k,2)=k(k+1)(2k+1)/6$. Higher order power sums
can be found recursively using \emph{Bernoulli numbers}.

The Bernoulli numbers (which are indexed by superscripts)
$1 = B^0$,  $B^1$, $B^2$, $B^3$, $B^4$,\dots{}
are defined by the formal equation
``$B^n = (B-1)^n$'' for  $n>1$,  where the
quotation marks indicate that the involved terms are to be
expanded in formal powers of $B$ before interpreting.
Thus:
\begin{itemize}
\item $B^2 = B^2 - 2B^1 + 1$, whence $B^1 = 1/2$,
\item $B^3 = B^3 - 3B^2 + 3B^1 - 1$, whence  $B^2 = 1/6$,
\end{itemize}
etc.
We thus get that $B^3 =0,$  $B^4 = -1/30$, $B^5 = 0$, $B^6 = 1/42$,
$B^7 = 0$, and so on.
It follows that for each $m$,
$$\p(k,m) = \frac{\mbox{\rm ``}(k+B)^{m+1} - B^{m+1}\mbox{\rm ''}}{m+1}$$
(\emph{Faulhaber's formula} \cite{CONGUY}.)
Thus, the coefficients $\p(k,m)$ can be efficiently calculated for small
values of $m$.
In particular, we get that
\begin{itemize}
\item $\p(k,3) = \frac{1}{4}k^4 + \half k^3 + \frac{1}{4}k^2$,
\item $\p(k,4) = \frac{1}{5}k^5 + \half k^4 + \frac{1}{3}k^3 - \frac{1}{30}k$,
\item $\p(k,5) = \frac{1}{6}k^6 + \half k^5 + \frac{5}{12}k^4 -
\frac{1}{12}k^2$,
\item $\p(k,6) = \frac{1}{7}k^7 + \half k^6 + \half k^5 -
\frac{1}{6}k^3 + \frac{1}{42}k$,
\item $\p(k,7) = \frac{1}{8}k^8 + \half k^7 + \frac{7}{12}k^6 -
\frac{7}{24}k^4 + \frac{1}{12}k^2$,
\end{itemize}
etc.
In order to show that this is enough, we need to bound the tail of the
series in
Equation \ref{sum}.
We will achieve this by effectively bounding the power sums.

\begin{lem}\label{upperbound}
Let $k$ be any natural number, and assume that
$f:(0,k)\to\R^+$ is such that $f''(x)$ exists, and is
nonnegative for all $x\in (0,k)$. Then:
$$\sum_{i=1}^{k} f(i) < \int_0^k f(x+\half) \dx.$$
\end{lem}
\begin{pf*}{Proof.}
For each interval $[i,i+1]$ ($i=0,\ldots,k-1$),
the tangent to the graph of $f(x+\half)$ at $x=i+\half$ goes
below the graph of $f(x+\half)$. This implies that the area of the
added part is greater than that of the uncovered part.
\qed\end{pf*}

Using Lemma \ref{upperbound}, we have that for all $m>1$,
$$\sum_{i=1}^{k-1} i^m < \int_0^{k-1} (x+\half)^m dx <
\frac{(k-\half)^{m+1}}{m+1}.$$
Thus,
\begin{eqnarray}
\sum_{m=N}^\infty \frac{\p(k-1,m)}{mn^m} < \sum_{m=N}^\infty
\frac{(k-\half)^{m+1}}{m(m+1)n^m} <
\frac{k-\half}{N(N+1)}\sum_{m=N}^\infty \(\frac{k-\half}{n}\)^m =
\nonumber\\
 = \frac{k-\half}{N(N+1)} \cdot
\frac{\(\frac{k-\half}{n}\)^N}{1-\frac{k-\half}{n}} =
 \frac{(k-\half)^{N+1}}{N(N+1)\(1-\frac{k-\half}{n}\)n^N}.
\end{eqnarray}

We thus have the following.

\begin{thm}\label{arbitrary}
Let $\pi^k_n$ denote the probability that all elements in a sample
of $k$ elements out of $n$
are distinct. For a natural number $N$, define
$$\epsilon^k_n(N):=
\frac{(k-\half)^{N+1}}{N(N+1)\(1-\frac{k-\half}{n}\)n^N}.$$
Then
$$\sum_{m=1}^{N-1} \frac{\p(k-1,m)}{mn^m}
<  -\ln(\pi^k_n) <
\sum_{m=1}^{N-1} \frac{\p(k-1,m)}{mn^m} + \epsilon^k_n(N).$$
\end{thm}

For example, for $N=2$ we get:
$$\frac{(k-1)k}{2n} <
-\ln(\pi^k_n) <
\frac{(k-1)k}{2n} + \frac{(k-\half)^3}{6n^2\(1-\frac{k-\half}{n}\)}.$$

We demonstrate the tightness of these bounds with a few concrete examples:

\begin{exmp}{\rm
Let us bound the probability that in a class of $5$ students there exist
two sharing the
same birthday. Using Theorem \ref{arbitrary} with $N=2$ we get by simple
calculation that
$\frac{2}{73} < -\ln(\pi^5_{365}) <
\frac{2}{73} + \frac{243}{2105320}$, or numerically,%
\footnote{
All calculations in this paper were performed using the GNU bc calculator
\cite{bc},
with a scale of $500$ digits.
}
$0.0273972 < -\ln(\pi^5_{365}) < 0.0275127$.
Thus, $0.0270253 < \beta^5_{365} < 0.0271377$.
Repeating the calculations with $N=3$ yields
$0.0271349 < \beta^5_{365} < 0.0271356$.
$N=4$ shows that $\beta^5_{365} = 0.0271355\ldots$.
}\end{exmp}

\begin{exmp}{\rm
We bound the probability that in a class of $73$ students there exist two
sharing the
same birthday, using $N=2$:
$\frac{36}{5} < -\ln(\pi^{73}_{365}) < \frac{36}{5}+\frac{121945}{255792}$,
and numerically we get that
$0.9992534 < \beta^{73}_{365} < 0.9995882$.
For $N=3$ we get
$0.9995365 < \beta^{73}_{365} < 0.9995631$,
and for $N=8$ we get that $\beta^{73}_{365} = 0.9995608\ldots$.
}\end{exmp}

In Theorem \ref{arbitrary}, $\epsilon^k_n(N)$ converges to $0$ exponentially
fast with $N$.
In fact, the upper bound is a very good approximation to the actual
probability, as can be seen
in the above examples. The reason for this is the effectiveness of the
bound in
Lemma \ref{upperbound} (see \cite{APPROX} for an analysis of this bound
as an approximation).

\bigskip

For $k<\sqrt{n}$, we can bound $\beta^k_n$ directly:
Note that for $|x|<1$ and odd $M$,
$\sum_{m=0}^M (-x)^m/m! < e^{-x} < \sum_{m=0}^{M+1} (-x)^m/m!$.
\begin{cor}\label{M}
Let $\beta^k_n$ denote the probability of a birthday surprise in a sample
of $k$ out of $n$,
and let $l_N(k,n)$ and $u_N(k,n)$ be the lower and upper bounds from
Theorem \ref{arbitrary},
respectively. Then for all odd $M$,
$$-\sum_{m=1}^{M+1} \frac{(-l_N(k,n))^m}{m!} < \beta^k_n < -\sum_{m=1}^{M}
\frac{(-u_N(k,n))^m}{m!}.$$
\end{cor}

For example, when $M=1$ we get that
\begin{equation}\label{M=1}
\frac{(k-1)k}{2n}-\frac{(k-1)^2k^2}{4n^2} < \beta^k_n < \frac{(k-1)k}{2n}+
\frac{(k-\half)^3}{6n^2\(1-\frac{k-\half}{n}\)}.
\end{equation}

The explicit bounds become more complicated when $M>1$,
but once the lower and upper bounds in Theorem \ref{arbitrary}
are computed numerically, bounding $\beta^k_n$ using Corollary \ref{M}
is easy. However, Corollary \ref{M} is not really needed in order to
deduce the bounds -- these can be calculated directly from the bounds
of Theorem \ref{arbitrary}, e.g.\ using the exponential function built in
calculators.

\begin{rem}\label{better}
{\rm
\begin{enumerate}
\item It can be proved directly that in fact $\beta^k_n < \frac{(k-1)k}{2n}$
\cite{BELLARE1}.
However, it is not clear how to extend the direct argument
to get tighter bounds in a straightforward manner.
\item Our lower bound in Equation \ref{M=1} compares favorably with the lower
bound
$\(1-\frac{1}{e}\)\frac{(k-1)k}{2n}$ from \cite{BELLARE1} when
$k\le \sqrt{2n/e}$
(when $k>\sqrt{2n/e}$ we need to take larger values of $M$ to get a
better approximation).
\item $\p(k-1,m)$ is bounded from below by $(k-1)^{m+1}/(m+1)$. This implies a
slight improvement on Theorem \ref{arbitrary}.
\end{enumerate}
}
\end{rem}

\section{Some applications}

\subsection{Arbitrary percision calculators}
\emph{Arbitrary percision calculators} do calculations
to any desired level of accuracy. Well-known examples are the \emph{bc}
and \emph{GNU bc} \cite{bc} calculators.
Theorem \ref{arbitrary} allows calculting $\beta^k_n$
to any desired level of accuracy (in this case, the parameter $N$
will be determined by the required level of accuracy), and in practical time.
An example of such calculation appears below (Example \ref{huge}).

\subsection{Cryptography}
The probability of a birthday surprise plays an important role in the security
analysis of various cryptographic systems. 
For this purpose,
it is common to use the approximation $\beta^k_n\approx k^2/2n$.
However, in \emph{concrete security} analysis it is prefered to have exact
bounds rather than
estimations (see \cite{BELLARE} and references therein).

The second item of Remark \ref{better} implies that security bounds derived
using
earlier methods are tighter than previously thought.
The following example demonstrates the tightness of the bounds of Theorem
\ref{arbitrary}
for these purposes.

\begin{exmp}\label{huge}
{\rm
In \cite{BELLARE2}, $\beta^{2^{32}}_{2^{128}}$ is estimated approximately.
Using Theorem \ref{arbitrary} with $N=2$, we get that in fact,
$$
2^{-65.0000000003359036150250796039103} <
\beta^{2^{32}}_{2^{128}} <
2^{-65.0000000003359036150250796039042}.$$
With $N=3$ we get that $\beta^{2^{32}}_{2^{128}}$ lies between
$$2^{-65.000000000335903615025079603904203942942489665995829764250752}$$
and
$$2^{-65.000000000335903615025079603904203942942489665995829764250713}.$$
The remarkable tightness of these bounds is due to the fact that
$2^{32}$ is much smaller than $2^{128}$.
}\end{exmp}

Another application of our results is for
estimations of the quality of approximations
such as ${n\choose k}\approx n^k/k!$ (when $k<<n$):
\begin{fact}
${n\choose k} = \frac{n^k}{k!}\cdot \pi^k_n$.
\end{fact}
Thus the quality of this approximation is directly related to
the quality of the approximation $\pi^k_n\approx 1$,
which is well understood via Theorem \ref{arbitrary}.

\bigskip

$\pi^k_n$ appears in many other natural contexts.
For example, assume that a function $f:\{0,\ldots,n-1\}\to \{0,\ldots,n-1\}$ is
chosen with uniform probability from the set of all such functions, and fix
an element $x\in \{0,\ldots,n-1\}$. Then we have the following immediate
observation.
\begin{fact}
The probability that the orbit of $x$
under $f$ has size exactly $k$ is $\pi^k_n \cdot \frac{k}{n}$.
The probability that the size of the orbit of $x$ is larger than
$k$ is simply $\pi^k_n$.
\end{fact}
These probabilities play an important role in the theory of iterative
pseudorandom number generation
(see \cite{ADI} for a typical example).

\section{Final remarks and acknowledgments}
For a nice account of power sums see \cite{EDWARDS}.
An accessible presentation and proof of Faulhaber's formula appears in
\cite{CONGUY}.
The author thanks John H.\ Conway for the nice introduction to Bernoulli
numbers,
and Ron Adin for reading this note and detecting some typos.

\end{document}